\documentclass[noamsfonts]{amsproc}
\usepackage{aspmproc}

\def\CC{{\bf C}}
\def\NN{{\bf N}}

\def\QQ{{\bf Q}}

\def\del{{\partial}}

\newcommand {\be} {\begin{equation}}
\newcommand {\eeqn} {\end{equation}}
\newcommand {\bea} { \begin{eqnarray}}
\newcommand {\eea} {\end{eqnarray}}
\newcommand {\beas} { \begin{eqnarray*}}
\newcommand {\eeas} {\end{eqnarray*}}
\newcommand {\ra} {\rightarrow}
\newcommand {\lra} {\longrightarrow}

\newtheorem {lemma} {LEMMA} [section]
\newtheorem {theorem}[lemma]{THEOREM}
\newtheorem {prop}[lemma]{PROPOSITION}

\newtheorem {definition}[lemma]{DEFINITION}
\newtheorem {question}[lemma]{QUESTION}

\title[A Fibration introduced via the Kobayashi pseudometric]
{The Kobayashi pseudometric on algebraic \\
manifold and a canonical fibration}

\author[S. S.Y. Lu]{Steven Shin-Yi Lu}

\address{
Department of Mathematics\\
the University of Toronto\\
100 St. George Street, 4th floor\\
Toronto, Ontario\\
Canada\\
$\ $\\
emails: slu@math.toronto.edu\\
$\ \ \ \ $ slu@mpim-bonn.mpg.de
}

\thanks{Partially supported by a Canadian NSERC grant, The Max Planck
Institute for Mathematics and Osaka University.}

\begin{document}

\begin{abstract}
Given a compact complex  manifold $X$ of dimension $n$,
we define a bimeromorphic invariant $\kappa'_+(X)$ as the maximum $p$ 
for which there is a saturated line subsheaf $L$ of 
the sheaf of holomorphic $p$~forms whose
Kodaira dimension $\kappa (L)$ equals $p$. We call $X$ special
if  $\kappa'_+(X)=0$. We observe from earlier works that among the
algebraic $X$ with $\kappa'_+(X)\neq 2$ the special ones
are in fact characterized by vanishing Kobayashi 
pseudometric (and holomorphic connectedness) for $n=2$ without 
any condition, for $n=3$ when $\kappa(X)\neq 0, 3$, for
birationally abelian fibrations over curves and for
rationally connected fibrations over any of these examples.
We show here that rationally connected fibrations over special
manifolds and algebraic manifolds
with zero Kodaira dimension are special. 

We use the well-known
construction of F. Campana to give, for each projective
$X$ a canonical fibration $f: X\rightarrow Y$ holomorphic outside
a proper subvariety of $Y$ and whose general fibers are special.
We show that the inherited orbifold structure on $Y$ defined via 
the minimum multiplicity of those of the components of each fiber
does not admit positive dimensional special sub-orbifolds through
the general points of $Y$.
We note that the Iitaka fibration or any rationally connected fibration
of $X$ is a natural factor of $f$ and we show that this solves a
general conjecture in Mori's classification program of algebraic 
varieties, namely,  that an algebraic variety
is either of general type, or (birationally)
has a canonical fibration with positive dimensional special 
type fibers that factors through the Iitaka and rationally 
connected fibrations of $X$. In fact, assuming the well-known
conjectures in Mori's Program but extended to log-varieties, which implies 
in particular that orbifolds of negative Kodaira dimension 
admits a fibration with positive dimensional special type fibers, 
our fibration is then nothing but the iterations of the orbifold
Iitaka and Mori's fibrations that end in a general type 
orbifold base. 

\end{abstract}

\maketitle

\section{Introduction and acknowledgements}

The purpose of this paper is twofold:  (1) to give a quick 
summary of our results on the Kobayashi pseudometric on algebraic
manifolds, especially on our characterization of those with
identically vanishing Kobayashi pseudometric that resulted
in the introduction of a natural birational invariant and its associated
canonical fibration, and (2) to show that this fibration resolves
the conjectured alternative in Mori's Program that an algebraic variety
$X$ has a canonical fibration with an orbifold base ``of general type'' 
(factoring through the Iitaka fibration or any rationally connected 
fibration of $X$ in the birational category) that has positive 
dimensional fibers if and only if $X$ is not of general type.\\

The author would like to thank J.P. Demailly, A. Fujiki,
B. Shiffman and M. Zaidenberg for their keen interest, initiation and 
guidance in this problem. A very special thanks is owed to  
F. Campana for his kindness in, among others, 
communicating his very stimulating preprint for which this
paper can only serve to add a few complementary results.
He also expresses his thanks to Professors Y. I. Manin, 
F. Hirzebruch and Dr. Y. Holla for some useful discussions 
at the Max-Planck-Institute in Bonn and finally to Y. Kawamata and 
F. Bogomolov for their suggestions and encouragements.

\section{Preliminaries on fibrations and orbifolds}

Since the key objects under study in this paper are all
natural bimeromorphic invariants, we  assume henceforth
that all objects are compact smooth manifolds and that every 
(meromorphic) map, unless otherwise specified,
is bimeromorphic to a flat morphism to the target space
and that the map is smooth (of maximal rank) outside a reduced
simple normal crossing divisor, the discriminant locus,
on the target manifold. We can do this
by 
the resolution of singularity theorem and the
flattening theorem \cite{Hi}.\\[-3mm]

Recall that sections, if exist, of powers of a 
holomorphic line bundle $L$ on a compact complex manifold $X$ 
gives rise to a fibration $$I_L : X \lra B=I_L(X),$$
possibly after replacing $X$ by a modification.
The Kodaira dimension $\kappa(L)$ of $L$ is defined to be $\dim I_L(X)$ 
in this case, and to be $-\infty$ otherwise. The general fibers $F$
of $I_L$ enjoys $\kappa(L|_F)=0$. 
These definitions generalize to torsion free line sheaves by taking
the double dual and to $\QQ$-line bundles, which are
line bundles tensored by $\QQ$-divisors. We set  $I_X:=I_{K_L}$  
and  call $\kappa(X):=\kappa(K_X)$ the Kodaira dimension of $X$.
The general fibers of $I_X$ have zero Kodaira dimension.

A (meromorphic) fibration $f: X\ra Y$ naturally imposes an orbifold 
structure on $Y$. To define it, we first replace 
$X$ by a modification to make $f$ holomorphic.
Given a reduced irreducible divisor $D$ on $Y$, we may
write $f^*D=\sum_i m_iD_i$ for $m_i\in \NN$ and reduced irreducible
divisors $D_i$ in $X$. Then we define the (minimum) multiplicity of 
$f$ over $D$ by
$$m^f(D)=\min \{\ m_i\ |\ f(D_i)=D\ \},$$
which is independent of the our choice of modification.
We note that the classical multiplicity $m_f(D)$ is defined
by replacing $\min$ above by $\gcd$ and is the more natural
multiplicity to consider for questions concerning the 
fundamental group, see \cite{Ca01}.
However, the multiplicity $m^f(D)$ is more natural for
our problems as will be indicated in this section. We define
$$D^f=\sum \{\ (1-\frac{1}{m^f(D)})D\ |\ 
D\ \mbox{a reduced irreducible divisor on}\ Y\ \}$$
and call the pair $(Y,D^f)$ an orbifold which we
also write as $(Y,f)$ or simply $Y^\del$. Note
that $D^f$ is supported on the discriminant locus of $f$.
If $C$ is a curve, we call a holomorphic map
$h: C\ra Y$ a  Z-orbifold map to $Y^\del$ if 
$g^*D\ge m^f(D)(g^*D)_{red}$ for all components $D$
in $D^f$. This gives us the notion of the orbifold
Kobayashi pseudometric $d_{Y^\del}$ (define!).

\begin{lemma}[\cite{Lu01}] With the above setup, assume that 
$d_F\equiv 0$ for the general fiber $F$ of $f$. Then every
point in $Y$ has a neighborhood $U$ with 
$d_{f^{-1}(U)}=f^*d_{U^\del}$. If moreover, $F$ is a
rationally connected fibration over a birationally abelian variety, 
then $d_X=f^*d_{Y^\del}$.
\end{lemma}

We define the orbifold canonical bundle of $Y^\del$ by
the $\QQ$-line bundle $K_{Y^\del}=K_Y(D^f)$,
understood as a  line bundle on Y after
a tensor power, and the orbifold Kodaira 
dimension $\kappa(Y^\del)$ of ${Y^\del}$ by $\kappa(K_{Y^\del})$.
If $\kappa(Y^\del)=\dim Y$, then $Y^\del$ and $f$ 
are said to be of general type. 
The usual lemma of Castelnuovo-De Franchis, 
as observed by Bogomolov (\cite{Bo}, \cite{Re}),
gives:

\begin{prop}[Castelnuovo-De Franchis, Bogomolov]\label{CD}
$\ \ \ \ \ \ $
Let $L$ be a saturated line subsheaf of $\Omega^p_X$
with $X$ projective. Then
\begin{itemize}
\item[(I)] $\kappa(L)\le p$
\item[(II)] If $\kappa(L)=p$, then $I_L$ defines a fibration from $X$
to a projective base $B$ of dimension $p$ and
$I_L^*K_B$ saturates to $L$ in $\Omega_X^p$. 
\item[(II')]
More  specifically, if $\kappa(L)=p$, 
then outside the singularity of the discriminant locus in B,
we have $L=I_L^*K_B(dI_L)$ where $dI_L$ is a  divisor given by
$\sum I_L^*(D_i) - I_L^*(D_i)_{red}$ and $dI_L\geq I_L^*D^{I_L}$ as
$\QQ$-divisors with equality achieved on every fiber. The latter
implies that $\kappa(B^\del)=\kappa(L)$ so that $B^\del$ is an
orbifold of general type.

\end{itemize}
\end{prop}

Parts (I) and (II) are found in \cite{Bo}, but see \cite{Re,Lu01}
for (II'). Note that a line sheaf
given by (II) and (II') gives rise to a fibration of general type
while the converse follows easily from our definition of $m^f$.

\begin{definition} We define $\kappa'_+(X)$ to be the maximum $p\geq 0$
for which a line subsheaf $L$ of $\Omega^p_X$ 
(the sheaf of holomorphic $p$-forms) exists with
$\kappa(L)=p$. $X$ is said to be special if $\kappa'_+(X)=0.$
\end{definition}

This invariant was introduced in \cite{Lu01} as a natural
refinement of an invariant $\kappa^+$ given in \cite{Ca}. Both
were really forced upon us by our respective problems. Indeed,
vanishing of $\kappa^+$ coincides precisely with the identical vanishing
of the Kobayashi metric in all the cases we are able to give
a characterization of the latter among projective varieties
and this remains true even in
the quasiprojective case where the logarithmic version of 
$\kappa^+$ is used.
The following is quoted from \cite{Lu01}.

\begin{theorem} Let $X$ be an $n$-dimensional projective manifold
not of general type.
We impose the condition  for $n=3$ that $\kappa(X)\neq 0$ and
$\kappa'_+(X)\neq 2$  while
for $n>3$, we let $X$ be birational to an abelian fibration 
over a curve or to a rationally connected fibration  
over any of the above examples. Then
$d_X\equiv 0$ if and only if $\kappa'_+(X)=0.$
\end{theorem}

\section{An orbifold Iitaka conjecture on Kodaira dimension}

We henceforth restrict to the projective category.
Fibrations whose general fibers are special are called
special fibrations. Two questions that arise naturally are whether
there is a unique line sheaf associated to $\kappa'_+(X)$
and whether the fibration(s) associated with $\kappa'_+(X)$
have special fibers. The answer to both questions has been 
checked to be in the affirmative for some classes of varieties, 
including those in the theorem above
(see \cite{Lu01}).
The two questions are indeed related
and their affirmative solution follows from the same for the following:

\begin{question}\label{conj1} Let $f:X\ra Y$, $g:Y\ra Z$ be fibrations
such that $h=g\circ f$ is a fibration of general type.
If the restrictions of $f$ to the general fibers of $h$ are also of
general type, is $f$ then of general type?
\end{question}

We address this question in \cite{Lu02}, where we prove
the orbifold versions of the theorems proved here.
In the case $f$ is the identity map and $Z$ is of general type,
theorems of Viehweg (see \cite{En}), Kawamata (\cite{Ka}) 
and Kollar (\cite{Ko}) all answered this question
in the affirmative.
These theorems are in turn special cases of the C$_{n,m}$ conjecture of 
Iitaka whose orbifold generalization asserts an affirmative 
answer to the above question.

We prove here the following case 
of this orbifold C$_{n,m}$ conjecture.

\begin{theorem}\label{Cnm} Let $f:X\ra Y$ and $g:Y\ra Z$ 
be projective fibrations with $g$ of general type.
If $\kappa(Y_z)\ge 0$ for the general fiber $Y_z$ of $g$, then 
$\kappa(Y^\del)= dim(Z)+\kappa(Y_z^\del)$,
where $Y_z^\del=(Y_z, f|_{Y_z})$.
\end{theorem}

This is essentially a direct corollary of the following
theorem, which generalizes slightly the results of the people
mentioned above. Details of this implication can be found
in \cite{Ca01}, which has the version of this theorem via
the gcd-defined orbifold base, (or \cite{Lu02}) as given in \cite{Ka, En}.

\begin{theorem}\label{wp} 
Assume that $f: X\ra Y$ is a fibration where $X, Y$ are 
projective. Set 
$K_{X/Y^\del}=K_X\otimes f^*K_{Y^\del}^{-1}$ as $\QQ$-line 
bundles. 
Let $m$ be a positive integer such that $mD^f$ is a Cartier divisor
on $Y$. Then $f_{*}(K_{X/Y^\del}^{\otimes m})$ is weakly 
positive in the sense of \cite{Vi}, defined also in \cite{En}.
\end{theorem}

Without the orbifold $\del$, this is the Kawamata-Viehweg's theorem.
Weak positivity is a generalization of the notion of nefness of
a vector sheaf over a curve and enjoys
the following properties:\\
(P1) A torsion free coherent sheaf is weakly positive if it 
contains a weakly positive subsheaf of the same rank.\\
(P2) Let $\Upsilon$ be a torsion free coherent sheaf on $Y$.
Let $v: Y'\ra Y$ be a flat base change. If $v^*\Upsilon$
is weakly positive, then so is $\Upsilon$.

{\bf Proof of Theorem~\ref{wp}}:
Recall that $D=D^f_{red}=\sum_i D_i$ is a 
simple normal crossing divisor on $Y$. 
Write $f^*D_i=\sum_j m_{ij}D_{ij}$ and let 
$M_i$ be the least common multiple of 
$\{\ m_{ij}\ |\ f(D_{ij})=D_i\ \}$. 
The Kawamata branched covering trick (\cite{Ka}) 
guarantees the existence of a finite Galois cover 
$v: Y'\ra Y$ with $Y'$ smooth such that $v^*D_i=M_iD'_i$ 
for some effective divisor. We let 
$u:X'\ra X$ be the composition $u=u_1\circ d$, where
   $u_1:X_1\ra X$ is the base change of $X$ by $v$
composed with the normalization map, and
$d:X'\ra X_1$ is a desingularization which 
is an isomorphism above the smooth locus of $X_1$.
We also denote by $f_1:X_1\ra Y'$ and $f'=f_1\circ d:X'\ra Y'$ the 
resulting pullbacks of $f$ by $v$. Note that if $f(D_{ij})=D_i$,
then $u^*(m_{ij}D_{ij})=M_iD'_{ij}$ 
for some effective divisor $D'_{ij}$ by construction.
Hence, applying $u^*$ to the $\QQ$-divisor $f^*D^f$ renders
it into an effective (integral) divisor on $X'$.

Now, a local computation shows the
following natural 
inclusion of torsion free coherent sheaves of
the same rank over $Y'$:  
$$
u_{1*}: (f')_{*}\big(K_{X'/Y'}^{\otimes m}\big)
\ra v^{*}f_{*}\big(K_{X/Y}(-f^*D^f)^{\otimes m}\big)
=v^{*}f_{*}\big(K_{X/Y^\del}^{\otimes m}\big).
\ \ \ \ \ (*)
$$
A round about way to see this is as follows:
This fact without the $(-f^*D^f)$ is well-known
(e.g., \cite{Ka}, \cite[Lemma 13]{En}). Hence, we
only need to check $(*)$ above an analytic neighborhood
$V$ in $Y$ of a generic point $p\in D'_i=v^{-1}D_i$. Let $D'_{ij}$
be a component of $f_1^{-1}D'_i$ dominating $D_{ij}$ and let
$U$ be a small neighborhood in $X_1$
of a generic point $q$ on $D'_{ij}$ above $p$. These neighborhoods
as $j$ varies cover $f_1^{-1}(V)$ outside
a codimension two subset if we include also the open subset
$f_1^{-1}(V\setminus D'_i)$, above which $K_{X_1/Y'}$ and 
$u_1^*K_{X/Y}$ are naturally identical. 
Since $d_*K_{X'/Y'}=K_{X_1/Y}$ on $X_1$,
see the proof of Lemma~13 in \cite{En}, and since $X_1$ is
normal, it is sufficient to show that $u_1^*K_{X/Y}(-u_1^*(df))$ is naturally
identified with $K_{X_1/Y'}$ on $U$, where 
$(df)=\sum f^*(D_i)-f^*(D_i)_{red}\geq f^*D^f$ is as given in 
part (II') of proposition~\ref{CD}. This identification is
easily seen from the fact that, on a neighborhood of $q$, we
can write $u=u_1=u_{ij}\circ u'_{ij}$ where $u_{ij}:X_{ij}\ra X$ is
obtained from the  base change by an $m_{ij}$-cyclic cover $Y_{ij}$ 
of $Y$ branched precisely on $D_i$ and $u'_{ij}$ is obtained
from the (remainder) cyclic cover of $Y'$ over $Y_{ij}$.
Let $f_{ij}:X_{ij}\ra Y_{ij}$ be the map corresponding to $f$.
Then by construction, $K_{X_{ij}}=u_{ij}^*K_X$ and $f_{ij}$ is
smooth on a neighborhood $W$ of $u'_{ij}(q)$. The first implies 
$K_{X_{ij}/Y_{ij}}=u_{ij}^*K_{X/Y}(-(df))$ on $W$ while the
second implies that $K_{X_1/Y'}=(u'_{ij})^*K_{X_{ij}/Y_{ij}}$ 
on $U\cap u_{ij}^{\prime -1}(W)$
as required for $(*)$. Combining the theorem of Kawamata-Viehweg
with $(*)$ via (P1) and (P2) now gives us theorem~\ref{wp}
and hence theorem~\ref{Cnm}.
$\qed$

By taking $X=Y$ in Theorem~\ref{Cnm}, we get the following 
corollaries.

\begin{theorem}\label{Cnm+} Let $g:Y\ra Z$ 
be a fibration of general type.
Then $\kappa(Y)\geq dim(Z)+\kappa(Y_z)$, where $Y_z$ is the 
general fiber of $g$. In particular, if 
$Y_z$ is of general type, then $Y$ is also.
\end{theorem}

\begin{theorem}\label{Sp} A projective manifold $X$ is special
if $\kappa(X)=0$.
\end{theorem}

{\bf REMARK:} We have been notified by Campana after 
this paper appeared as an MPI preprint that he, in the latest 
versions of his preprint cited here, has also adopted the 
``min'' definition for the orbifold base and strengthened 
his theorems to the same ones given in this section.

\section{A canonical special fibration}
Due to page limitations, this section will be brief. 
%

We first observe that the general fiber of any fibration
has the same $\kappa'_+$ as one can deduce from 
semi-continuity and other facts in the subject (\cite{Lu02}).
Let $X$ be a smooth projective variety and $T$
the union of all irreducible components of 
the Chow scheme $Chow(X)$ whose general points $t$
correspond to irreducible reduced special 
subvarieties $V_t$ of $X$, including the trivial ones.
$T$ gives an equivalence relation on $X$ where
$x,x'\in X$ are $T$-equivalent if and only if they  are
contained in a $T$-chain (i.e., a connected union of
finitely many of the $V_t$'s). We say that $X$ is $T$-connected
if all pairs of points in $X$ are $T$-equivalent. If $X$ is
$T$-connected, then $X$ is special: For if $f:X\ra B$ is a
fibration of general type, then $T$-connectedness gives us
a family $Y\ra S\ra T$ such that the general $Y_s$ identifies 
with a special subvariety in $X$ and that the composition 
$Y\ra X\ra B$ is generically finite, showing that $Y$
is of general type and hence so are the general fibers
$Y_s$, a contradiction. Hence, the construction
of Campana \cite[\S 14]{Ca01, Ca04}, whose notation we also
adopt here as ambiguity no longer exists, 
see also \cite{KoMM,Kol2,De}, gives:

\begin{theorem}\label{Main} There is a unique fibration
$c_X: X\ra C(X)$, holomorphic outside a normal crossing
divisor in $C(X)$, such that

{\rm (1)} The general fiber of $c_X$ is special.

{\rm (2)} If $Z\subset X$ is a special subvariety of $X$ going
through a general point $a\in X$, then $Z$ is contained in the fiber
   of $c_X$ through $a$.

{\rm (3)} If $X$ is not of general type, then $c_X$ has positive
fiber dimension.

{\rm (4)} The general point of the orbifold base $C(X)^\del$ 
is not contained in a nontrivial special orbifold subvariety
(definition given in the proof).
\end{theorem}

{\bf Proof}: The first parts can be found in the references
above, so we only need to show parts (3) and (4) here. By the assumption 
of (3), $\kappa'_+(X)<\dim X$ so that we have a fibration $f:X\ra Z$
of general type with positive fiber dimension. Now the general
fibers $X_z$ of $X$ are not of general type, for otherwise $X$ 
would be of general type by theorem~\ref{Cnm+}. Hence we can
repeat this process with $X_z$ until we end with a positive 
dimensional special subvariety through the general point of $X$.
For (4), we define $\kappa'_+$ for an orbifold $X^\del$ by the
same definition as for $\kappa'_+(X)$ except we replace 
subsheaves of $\Omega_X$ with subsheaves of 
a canonical locally free subsheaf $\Omega_{X^\del}$ of 
$\Omega_{X'}$ ($X'$ being a Galois cover that kills the orbifold
structure of $X$) invariant under the Galois group.  We call 
X special if $\kappa'_+(X^\del)=0$. Then (4) follows from 
the following easily deduced 
lemma. See \cite{Lu02} for more details on these.

\begin{lemma}\label{s/s} $X$ is special if it has a special fibration 
$f: X\ra Y$ such that  $Y^\del$, understood as $(Y,D^f)$, is a 
special orbifold. \qed
\end{lemma}

Now, a projective variety $X$ that is rationally connected 
(having a rational curve through every pair of points of $X$) 
or has zero Kodaira dimension is
special by above and theorem~\ref{Sp}, respectively. Also,
a rationally connected fibration is a fibration whose 
general fibers are rationally connected.
Hence (1) and (2) of theorem~\ref{Main} shows that the Iitaka 
and rationally connected fibrations of $X$ are factors of
$c_X$. This coupled with (3) solves the Mori's alternative 
stated at the introduction. 

Also,  assuming the standard conjectures in Mori's program, 
which asserts in particular that orbifolds of 
negative Kodaira dimension admits a ``rationally connected'' 
fibration (called Mori fibration)
with positive dimensional fibers,
then lemma~\ref{s/s} allows us to deduce that $c_X$ is
just the iteration of the orbifold Iitaka and Mori
fibrations that ends at a general type orbifold base, see
\cite{Lu02} for the details.\\



\end{document}